\documentclass[a4paper,oneside,11pt]{article}

\usepackage{amsfonts,color}

\newtheorem{dfn}{Definition}[section]
\newtheorem{tw}[dfn]{Theorem}
\newtheorem{prop}[dfn]{Proposition}
\newtheorem{rem}[dfn]{Remark}

\newtheorem{lem}[dfn]{Lemma}

\usepackage{amssymb} 
\usepackage{amsmath}
\numberwithin{equation}{section}

\usepackage{anysize}
\usepackage{array}
\usepackage{enumerate}

\author{Micha\l \ Barski  \\ \small  Faculty of Mathematics, Warsaw University, Poland\\
 \small{\it m.barski@mimuw.edu.pl} \bigskip \\
\\
Jerzy Zabczyk
\\ \small Institute of Mathematics, Polish Academy of
Sciences,
     Warsaw,  Poland\\ \small{\it zabczyk@impan.pl}}

\title{\bf On  generalized CIR equations}

\begin{document}

\maketitle

\begin{abstract}The paper is concerned with  stochastic equations for the short rate process $R$ 
$$
dR(t)=F(R(t))dt+G(R(t-))dZ(t),
$$
in the affine model of the bond prices. The  equation is driven by a L\'evy martingale $Z$.  It is shown that  the discounted bond prices are local martingales if  either  $Z$ 
is a stable process of index $\alpha\in(1,2]$,\,$F(x)= ax +b,\, b\geq 0$,\,\,  $G(x)=cx^{1/\alpha},\,c>0$ or $Z$ must be  a L\'evy martingale with positive jumps and trajectories of bounded variation, $F(x)= ax +b,\, b\geq 0$\,\,
and  G is a constant. The result generalizes
the well known Cox-Ingersoll-Ross result from \cite{CIR} and extends the Vasi\v cek result,  see \cite{Vasicek}, to non-negative short rates.

\end{abstract}

\noindent
\begin{quote}
\noindent \textbf{Key words}: CIR model, bond market, HJM condition, stable martingales.

\textbf{AMS Subject Classification}: 60G60, 60H20, 91B24 91B70.

\textbf{JEL Classification Numbers}: G10,G12.
\end{quote}

\bigskip
\section{Introduction}
The Cox-Ingersoll-Ross (CIR) equation is a stochastic equation for a non-negative process $R(t), t\geq 0$, of the form
\begin{gather}\label{rownanie na R}
dR(t)=(aR(t)+b)dt+c\sqrt{R(t)}dW(t), \quad R(0)=x\geq 0,
\end{gather}
where $a,b,c$ are constants, $b\geq 0$, $c\geq 0$, and $W$ is a Wiener process. It was introduced in the paper \cite{CIR} to model short
rate process $R(t)$ for which the bond prices $P(t,T), 0\leq t\leq T<+\infty$ have an affine structure, that is
\begin{gather}\label{rownanie na P}
P(t,T)=e^{-A(T-t)-B(T-t)R(t)}, \quad 0\leq t\leq T,
\end{gather}
where $A$ and $B$ are smooth, non-negative, deterministic functions and the discounted bond prices are  local martingales. That is the following   {\it martingale property}   is satisfied:
\begin{align*}
\phantom{(MP)}\hskip2ex& \text{For each $T>0$,\, the discounted bond price process}\\
\bf{(MP)}\hskip4ex& \hskip20ex\hat{P}(t,T)=P(t,T) e^{-\int_0^t R(s))ds},\quad\ t\in[0,T],\\
\phantom{(MP)}\hskip2ex& \text{is a local martingale}.
\end{align*}
This property implies that
the market with bonds priced by \eqref{rownanie na P} does not allow arbitrage. In fact a stronger result is true,  see e.g. Filipovi\'c \cite{FilipovicATS}.
If one requires that the short rate process $R$ is only a time homogeneous path continuous  Markov process and the bond prices
\eqref{rownanie na P}, with some functions $A$ and $B$, satisfy the martingale property, then it can be identified with the solution
of the CIR equation.

In the present paper we want to find all equations of the form
\begin{gather}\label{rownanie ogolne}
dR(t)=F(R(t))dt+G(R(t-))dZ(t), \quad R(0)=x\geq 0,
\end{gather}
where $Z$ is a real L\'evy martingale, such that the solution $R$ is non-negative and determines affine bond prices with the martingale property (MP).
Although there exists a precise description of the infinitesimal generator of such process due to Filipovi\'c, see \cite{FilipovicATS},
the analog of the CIR equation, in a more general setting of c\`adl\`ag short rate process $R$, is appearing in the present paper
for the first time. Our main theorems  give a precise characterization of all those  equations \eqref{rownanie ogolne} for which  the martingale property holds.  In our approach we do not use  the Filipovi\'c characterization but examine the problem via the Heat-Jarrow-Morton (HJM)
conditions for the absence of arbitrage. In the book in preparation \cite{BarskiZabczyk} some results using  the Filipovi\'c characterization are presented, but are  less complete. In \cite{BarskiZabczyk} one can also find a discrete time version of the problem in the framework of Filpovi\'c and Zabczyk, see \cite{FilipovicZabczyk}.

The choice of a L\'evy process in the equation  \eqref{rownanie ogolne} is motivated by the requirement that solutions to
\eqref{rownanie ogolne} should be  Markov processes. It is, however,  natural, in the light of a meta-theorem, that  Markov processes can be represented as  solutions
to stochastic equation with the noise  being L\'evy process but, in general,  infinite dimensional, see the book \cite{PeszatZabczyk}.\

In the present paper we allow only a one dimensional noise process. For the multidimensional L\'evy noise $Z = (Z_1 ,Z_2 ,\ldots ,Z_d)$ the problem of determining all equations
\begin{equation}
dR(t) = F(R(t))dt + \sum_{j=1}^{d} G_{j}(R(s-))d Z_{j}(s),\,\,\,\,R(0)=x\geq 0,
\end{equation}
for which the martingale property holds, is open. We  have  only partial answers which can be found in \cite{BarskiZabczyk}. It seems that even equations with multidimensional noise are not enough to represent all Markovian short rate processes  with the martingale property. 

As a part of the main proofs, in the section devoted to auxiliary results, we derive rather general results on positive invariance of stochastic equations and on equations with the martingale property, see Proposition \ref{positivity1}, Proposition \ref{tw G(0)=0}, Proposition \ref{tw char. ATS przyp ogolny}. They might be helpful for an  attempt to cover the multidimensional case.

\section{Results}

Our aim is to determine functions $F,G$ and a L\'evy process $Z$ in \eqref{rownanie ogolne}
as well as the functions $A,B$ in \eqref{rownanie na P} such that
solutions $R$ to the equation \eqref{rownanie ogolne} starting from arbitrary $x\geq 0$ are non-negative and the processes
$\hat{P}(t,T), t\in[0,T]$ are local martingales for arbitrary $T>0$. If this is the case then we say that
the equation \eqref{rownanie ogolne} determines an affine arbitrage free market or equivalently that  {\it equation \eqref{rownanie ogolne} has the  martingale property (MP)}.
\vskip2mm

\noindent It is natural to assume that
$$
0< P(t,T)\leq 1, \quad P(T,T)=1, \quad 0\leq t\leq T,
$$
and that $P(t,T)$ is a decreasing function of $T$. Since $R(t), t\geq 0$, should be a non-negative process,
$A$ and $B$ should be increasing, non-negative functions starting from zero, i.e.
\begin{gather}\label{zerowanie w zerze A i B}
A(0)=0, \quad B(0)=0.
\end{gather}
In the sequel we will argue that necessarily
\begin{gather}\label{pochodne w zerze A i B}
A^\prime(0)=0, \quad B^\prime(0)=1,
\end{gather}
see the reasoning preceding formula \eqref{warunki na pochodne A i B w zerze} in the sequel. Without loss of generality we will look for functions $G$ such that
\begin{gather}\label{warunek na G(bar x)}
\exists \bar{x}>0  \quad G(\bar{x})>0.
\end{gather}
We can reduce the general consideration to this case by multiplying
$G$ and $Z$ by $-1$, if $G$ is always non-positive. We also assume that
the L\'evy process is integrable with mean zero, i.e.
\begin{gather}\label{Z - intgrability}
\mathbb{E}\mid Z(t)\mid<+\infty, \quad \mathbb{E}Z(t)=0, \quad t\geq 0.
\end{gather}
Clearly, \eqref{Z - intgrability} means that $Z$ is a martingale. Note that if
$\mathbb{E}Z(t)\neq 0$, then $\mathbb{E}Z(t)=at, t\geq 0$ for some $a\neq 0$, and replacing
$Z$ by $\bar{Z}(t)=Z(t)-at, t\geq 0$, one arrives at the, equivalent to \eqref{rownanie ogolne}, equation
\begin{gather*}
dR(t)=(F(R(t))+a G(R(t)))dt+G(R(t-))d\bar{Z}(t), \quad R(0)=x\geq 0.
\end{gather*}
By $J$ we denote the Laplace exponent of the process $Z$, that is
\begin{gather}\label{Laplace exponent}
\mathbb{E}[e^{-\lambda Z(t)}]=e^{tJ(\lambda)}, \quad t\geq 0, \lambda\in\Lambda,
\end{gather}
where $\Lambda$ is the set for which the left side of \eqref{Laplace exponent} is well defined. It is well known that
\begin{gather}\label{J wzor}
J(\lambda)=-a\lambda+\frac{q}{2}\lambda^2+\int_{-\infty}^{+\infty}(e^{-\lambda y}-1+\lambda y\mathbf{1}_{(-1,1)}(y))\nu(dy),
\end{gather}
where $a\in\mathbb{R}$, $q\geq 0$ are constants and $\nu$ is a measure on $\mathbb{R}\setminus \{0\}$ such that
$$
\int_{-\infty}^{+\infty}(y^2\wedge 1)\nu(dy)<+\infty,
$$
and that
$$
\Lambda=\{\lambda \in\mathbb{R}: \int_{\mid y\mid\geq 1}e^{-\lambda y}\nu(dy)<+\infty\}.
$$
The measure $\nu$ is called the L\'evy measure of the process $Z$ and this measure together with the constants $a,q$
in \eqref{J wzor} determine the process $Z$ in a unique way. The fact that $Z$ is a martingale is equivalent to:
\begin{gather}\label{Z mart equivalent conditions}
\int_{\mid y\mid>1}\mid y\mid\nu(dy)<+\infty, \quad a+\int_{\mid y\mid>1}y\nu(dy)=0.
\end{gather}
If the L\'evy measure
of $Z$ is given by
$$
\nu(dy)= \frac{1}{y^{1+\alpha}}\mathbf{1}_{(0,+\infty)}(y)dy,
$$
then $J(\lambda)=c_{\alpha}\lambda^{\alpha}, z\geq 0$, where $\alpha\in(1,2)$ and $c_\alpha:=\frac{1}{\alpha(\alpha-1)}\Gamma(2-\alpha)>0$.
Here $\Gamma$ stands for the Gamma function. Then $Z$ is a stable martingale with index $\alpha$, providing that $a$ is such that \eqref{Z mart equivalent conditions} holds. This martingale will be denoted by $Z^\alpha$ for $\alpha\in(1,2)$ and $Z^\alpha$ with $\alpha=2$ stands for the Wiener process.

\begin{tw}\label{tw glowne necessity}
[Necessity]

\noindent Assume that the equation \eqref{rownanie ogolne} with functions $F,G$ which are continuous on $[0,+\infty)$
has the martingale property (MP) with some functions $A,B$ satisfying \eqref{zerowanie w zerze A i B} and \eqref{pochodne w zerze A i B}.
\vskip2mm
\noindent
\begin{enumerate}[I)]
\item If $G$ is differentiable on $(0,+\infty)$ and $G(\bar{x})>0$, $G^\prime(\bar{x})\neq 0$ for some
$\bar{x}>0$, then
\begin{enumerate}[a)]
\item $Z=Z^\alpha$ is a stable L\'evy process with index $\alpha\in(1,2]$ with positive jumps only,
\item $F(x)=ax+b$ with $a\in\mathbb{R}$, $b\geq 0$, $x\geq 0$,
\item $G(x)=c^\frac{1}{\alpha}x^{\frac{1}{\alpha}}$, $c>0$, $x\geq 0$.
\end{enumerate}
\item If G is a positive constant $\sigma$,  then
\begin{enumerate}[d)]
\item $Z$ has no Wiener part, i.e. $q$ in \eqref{J wzor} disappears,
\end{enumerate}
\begin{enumerate}[e)]
\item the martingale Z has positive jumps only and $\int_{0}^{+\infty} y\nu(dy) < + \infty$,
\end{enumerate}
\begin{enumerate}[f)]
\item $F(x)=ax+b$,\,\, $x\geq 0$, with $a\in\mathbb{R}$, $b\geq \sigma \int_{0}^{+\infty} y\nu(dy)$.
\end{enumerate}
\end{enumerate}
\end{tw}

Note that if there exists a number  $\bar{x}>0$ such that  $G(\bar{x})>0$, but for all such numbers   $G^\prime(\bar{x})= 0$, then  the function $G$ should be a positive constant on $[0, +\infty )$. This case is covered by Part II) of Theorem \ref{tw glowne necessity}.

\begin{tw}\label{tw sufficiency}
[Sufficiency]
\begin{enumerate}[I)]
\item The equation
\begin{gather}\label{rownanie 1002}
dR(t)=(aR(t)+b)dt+c^{\frac{1}{\alpha}}R(t-)^{\frac{1}{\alpha}}dZ^\alpha(t), \quad R(0)=x\geq 0, \quad a\in\mathbb{R}, b\geq 0, \ c>0,
\end{gather}
has a unique non-negative strong solution and satisfies the martingale property (MP).
The functions $A,B$ in \eqref{rownanie na P} are such that $B$ solves the equation
\begin{gather}\label{eq. for B}
B^{\prime}(v)=-cc_{\alpha}B^\alpha(v)+aB(v)+1, \quad v\geq 0,\quad B(0)=0,
\end{gather}
and $A$ is given by $A^{\prime}(v)=bB(v), v\geq 0,\,A(0)=0$.
\item If G is a positive constant $\sigma$ and $(d),\,(e), \ (f)$ in Theorem \ref{tw glowne necessity}, hold,  then  the equation
$$
dR(t) = (aR(t) + b) + \sigma dZ(t),\,\,\,R(0)=x \geq 0,\,\,t>0,
$$
has the martingale property (MP) and its solutions are non-negative processes. Moreover, $A$, $B$ are given by
\begin{align}
\label{eq. for B1} B'(v) =& B(v)a +1 ,\,\,B(0)=0,  \\
A'(v) =& B(v) (b- \sigma \int_0^{+\infty} y\nu(dy)) + \int_0^{+\infty}(1- e^{-\sigma B(v)y} )\nu(dy),\,\, A(0)=0.
\end{align}
\end{enumerate}
\end{tw}

\noindent The equation \eqref{eq. for B},  can be solved explicitelly for several values of  $\alpha$.   If $\alpha=2$ the equation   \eqref{rownanie 1002} becomes the CIR equation and
\eqref{eq. for B} boils down to the Riccati equation.

\begin{rem}
It is an open problem to characterize equations
\begin{equation}
dR(t)= F(R(t))dt + <G(R(s-)), dZ(s)>, \,\,R(0)=x \geq 0,
\end{equation}
where $Z$ is a L\'evy martingale of a given, greater than $1$, dimension, which solutions determine affine term structure with the martingale property.

It is easy to check that if $\alpha,\,\beta \in (1,2]$,\,\, $Z^\alpha,\,\, Z^\beta$ are independent stable martingales with indexes $\alpha,\,\,\beta$ respectively
and if $b\geq 0$ then the equation
\begin{gather}\label{Multi1}
dR(t)=(aR(t) + b)dt+R(t-)^{\frac{1}{\alpha}}dZ^\alpha(t) + R(t-)^{\frac{1}{\beta}}dZ^\beta(t),\,\,R(0)=x \geq 0,
\end{gather}
has a unique non-negative solution which determines an affine term structure with the martingale property.
\end{rem}

\section{Auxiliary results}
Here we derive some auxiliary results used in the proof of Theorem \ref{tw glowne necessity}.
\vskip2mm

\noindent Short rates should be non-negative processes. It is not surprising that in  the proof of necessary conditions in Theorem  
\ref{tw glowne necessity} we will need some results on equations which solutions are   positive invariant. The positive invariance has consequences for the structure of the L\'evy process $Z$, see  Proposition \ref{positivity1} and on the diffusion coefficient $G$, see   
Proposition \ref{tw G(0)=0}. We also derive analytic conditons on the coefficients of the  equation \eqref{rownanie ogolne} implied by the martingale property, see Proposition \ref{tw char. ATS przyp ogolny}.

\subsection{Positive invariance and  the noise process}

We show that if all solutions of  \eqref{rownanie ogolne} starting from non-negative points are non-negative then
 the noise process $Z$ can't have arbitrarily large negative jumps. This in turn affects properties of the Laplace exponent of $Z$.

 Below $supp(\nu)$ stands for the closed support of $\nu$.
\begin{prop}\label{positivity1}
Assume that $G(\bar{x})>0$ for some $\bar{x}>0$ and solutions to \eqref{rownanie ogolne} with all possible initial values
$x\geq 0$ are non-negative. Then
\begin{enumerate}[a)]
\item the support of $\nu$ satisfies
\begin{gather}\label{no large negative jumps}
supp(\nu)\subseteq (b,+\infty), \quad \text{for some} \quad b\in\mathbb{R},
\end{gather}
\item the domain $\Lambda$ of the Laplace exponent of $Z$ satisfies
$$
\Lambda\supseteq [0,+\infty),
$$
\item $J^\prime(\lambda)$ is finite for $\lambda>0$ and
$$
J^\prime(0+)=0.
$$
\end{enumerate}
\end{prop}
\noindent
{\bf Proof:} $a)$ For arbitrary negative $a$ let us consider the decomposition
$$
\nu(dy)=\nu_1(dy)+\nu_2(dy), \quad \nu_2(dy)=\mathbf{1}_{(-\infty,a)}(y)\nu(dy).
$$
The process $Z$ can be represented as a sum of two independent L\'evy processes $Z_1$, $Z_2$
such that $Z_2$ is a compound Poisson process with L\'evy measure $\nu_2$. The process $Z_2$
is a pure jump process with all jumps smaller than $a$. Let $X_1$ be a solution to the equation
$$
dX_1(t)=F(X_1(t))dt+G(X_1(t-))dZ_1(t), \quad X_1(0)=\bar{x},
$$
and let $c, r$ be positive numbers such that
$$
G(x)\geq c, \quad \text{for} \ x\in I:=\{z\geq 0: \bar{x}-r<z<\bar{x}+r\}.
$$
Let $\tau_1$ and $\tau_2$ be, respectively, the first exit time of $X_1$ from the interval
$I$ and the moment of the first jump of the process $Z_2$. There exists $T>0$ such that
$$
\mathbb{P}(\tau_1>T)>0, \quad \mathbb{P}(\tau_2<T)>0.
$$
Since the events $\{\tau_1>T\}$, $\{\tau_2<T\}$ are independent,
$$
\mathbb{P}(\{\tau_1>T\}\cap \{\tau_2<T\})=\mathbb{P}(\tau_1>T)\cdot \mathbb{P}(\tau_2<T)>0.
$$
Note that the solution of
\begin{gather*}
dX(t)=F(X(t))dt+G(X(t-))d(Z_1(t)+Z_2(t)), \quad X(0)=\bar{x},
\end{gather*}
satisfies
$$
X(t)=X_1(t), \quad \text{for} \ t<\tau_2.
$$
However, on the set $\{\tau_1>T\}\cap \{\tau_2<T\}$ we have
\begin{align*}
\triangle X(\tau_2)&=G(X(\tau_2-))\triangle Z(\tau_2)\\[1ex]
&=G(X(\tau_2-))\triangle Z_2(\tau_2),
\end{align*}
and
\begin{align*}
X(\tau_2)&=X_1(\tau_2-)+G(X_1(\tau_2-))\triangle Z_2(\tau_2)\\[1ex]
&\leq \bar{x}+r+ca.
\end{align*}
Thus, since $a$ can be chosen arbitrarily, we arrive to the condition $X(\tau_2)<0$, which is a contradiction.

\noindent Proofs of\,\,$b)$ and $c)$. Since the jumps of $Z$ are bounded from below by $b$, the integral
$$
\int_{\mid y\mid>1}e^{-\lambda y}\nu(dy)=\int_{\mid y\mid>1, y>b}e^{-\lambda y}\nu(dy),
$$
is finite for $\lambda\geq 0$. Consequently, the Laplace exponent $J(\lambda)$ is well defined for $\lambda\geq 0$.
Moreover, by \eqref{Z mart equivalent conditions} also
$$
J^\prime(\lambda)=-a+q\lambda+\int_{b}^{+\infty}y(\mathbf{1}_{\{\mid y\mid>1\}}-e^{-\lambda y})\nu(dy),
$$
is finite for $\lambda >0$ and
$$
\lim_{\lambda \downarrow 0} {\frac{J(\lambda)}{\lambda}} =J^\prime(0+)=0.
$$

\hfill$\square$

\subsection{Positive invariance and the diffusion coefficient}

We concentrate now on the equation
\begin{gather}\label{rownanie ogolne alfa stable}
dR(t)=F(R(t))dt+G(R(t-))dZ^\alpha(t), \quad R(0)=x\geq 0,
\end{gather}
where $Z^\alpha$ is a stable martingale with index $\alpha\in(1,2)$ and positive jumps only. Recall, its L\'evy measure has the form
$$
\nu(dy)=\mathbf{1}_{(0,+\infty)}\frac{1}{y^{1+\alpha}} dy.
$$
We will prove the following result:
\begin{prop}\label{tw G(0)=0}
	If  $G$ is a Lipschitz function and the equation \eqref{rownanie ogolne alfa stable}
	\,is positive invariant, then $G(0)=0$. 
\end{prop}
In the proof we use the classical maximal inequality
\begin{gather}\label{max inequality}
\mathbb{P}(\sup_{s\in[0,t]}\mid X(s)\mid\geq r)\leq \frac{3}{r}\mathbb{E}\mid X(t)\mid, \quad t>0,
\end{gather}
where $X$ is a c\`adl\`ag submartingale, see Proposition 7.12 in \cite{Kallenberg}.

\noindent We will also need the following lemma:
\begin{lem}\label{Lemat}
	For $2\geq p>\alpha>1$
	\begin{gather}\label{nier izomt}
	\mathbb{E}\mid\int_{0}^{t} g(s)dZ^\alpha_0(s)\mid^p\leq \frac{c_p}{p-\alpha}\mathbb{E}\int_{0}^{t}\mid g(s)\mid^p ds, \quad t\geq 0,
	\end{gather}
	with some $c_p>0$.
\end{lem}
\noindent
Here $Z^{\alpha}_0$ is a modified $\alpha$-stable martingale $Z^{\alpha}_0$ with L\'evy measure
$$
\nu(dy)=\mathbf{1}_{(0,1)}\frac{1}{y^{1+\alpha}} dy.
$$
Its jumps are thus bounded by $1$ and it is  identical with the process $Z^\alpha$ on the interval $[0,\tau_1)$, where $\tau_1$ is the first jump of $Z^\alpha$ exceeding $1$.

\vskip2mm
\noindent {\bf Proof of Lemma \ref{Lemat}:} Since the quadratic variation of the integral  $\int g(s)dZ^\alpha_0(s)$ equals
$$
\Big[\int g(s)dZ^\alpha_0(s)\Big](t)
=\int_{0}^{t}\int_{0}^{1} g^2(s)y^2\pi_0(ds,dy)
$$
where $\pi_0$ stands for the jump measure of $Z^\alpha_0$, by the Burkholder-Davis-Gundy inequality we obtain, for some $c_p>0$,
\begin{align*}
\mathbb{E}\mid \int_{0}^{t}g(s)dZ^\alpha_0(s)\mid^p&\leq c_p \mathbb{E}\Big[\int g(s)dZ^\alpha_0(s)\Big]^{\frac{p}{2}}(t)\\[1ex]
&=c_p\mathbb{E}\Big(\int_{0}^{t}\int_{0}^{1} g^2(s)y^2\pi_0(ds,dy)\Big)^{\frac{p}{2}},
\end{align*}
and further, since $p/2\leq 1$,
\begin{align*}
\mathbb{E}\mid \int_{0}^{t}g(s)dZ^\alpha_0(s)\mid^p
&\leq c_p \mathbb{E}\int_{0}^{t}\int_{0}^{1}\mid g(s)\mid^p y^p ds\frac{1}{y^{1+\alpha}}dz \\[1ex]
&\leq c_p\mathbb{E}\int_{0}^{t}\mid g(s)\mid^p ds\cdot \int_{0}^{1}\frac{y^p}{y^{1+\alpha}}dy\\[1ex]
&\leq c_p\mathbb{E}\int_{0}^{t}\mid g(s)\mid^p ds\cdot \int_{0}^{1}\frac{y^p}{y^{1+\alpha-p}}dy\\[1ex]
&\leq \frac{c_p}{p-\alpha}\mathbb{E}\int_{0}^{t}\mid g(s)\mid^p ds.
\end{align*}
\hfill $\square$

\vskip2mm
\noindent {\bf Proof of Proposition \ref{tw G(0)=0}:} We adopt the proof of Milian \cite{Milian} for the Wiener noise, which goes back to Gihman, Skorohod \cite{GihmanSkorohod}. Let us consider \eqref{rownanie ogolne alfa stable} with $x=0$. Then we can write $R$ in the form
$$
R(t)=\int_{0}^{t}F(R(s))ds+\int_{0}^{t}G(R(s-)-G(0))dZ^\alpha(s)+G(0)Z^\alpha(t), \quad t>0.
$$
Dividing by $t^{\frac{1}{\alpha}}$ yields
\begin{gather}\label{rownanie do analizy}
\frac{1}{t^{\frac{1}{\alpha}}}R(t)=\frac{1}{t^{\frac{1}{\alpha}}}\int_{0}^{t}F(R(s))ds+\frac{1}{t^{\frac{1}{\alpha}}}\int_{0}^{t}G(R(s-)-G(0))dZ^\alpha(s)+\frac{1}{t^{\frac{1}{\alpha}}}G(0)Z^\alpha(t), \quad t>0.
\end{gather}
Since
$$
\underset{t\rightarrow 0}{\lim\inf} \frac{1}{t^{\frac{1}{\alpha}}}Z^\alpha(t)=-\infty, \quad \underset{t\rightarrow 0}{\lim\sup} \frac{1}{t^{\frac{1}{\alpha}}}Z^\alpha(t)=+\infty,
$$
see \cite{Bertoin}, Theorem 5 in Section VIII, the last term in \eqref{rownanie do analizy} becomes negative for some sequence $t_n\downarrow 0$ providing that $G(0)\neq 0$. Since
$$
\frac{1}{t^{\frac{1}{\alpha}}}\int_{0}^{t}F(R(s))ds\underset{t\rightarrow 0}{\longrightarrow} 0,
$$
the assertion is true if we show that
$$
\frac{1}{t^{\frac{1}{\alpha}}}\int_{0}^{t}G(R(s-)-G(0))dZ^\alpha(s)\underset{t\rightarrow 0}{\longrightarrow} 0.
$$
Let us denote $g(s):=G(R(s-)-G(0)$. In the neighborhood of zero we can replace $Z^\alpha$ by
$Z^{\alpha}_0$. Then, by \eqref{max inequality}, for the submartingale $\mid\int_{0}^{s}g(u)dZ^\alpha_0(u)\mid^p$, with $2>p>\alpha>1$, we have
\begin{align}\label{nier. 1}\nonumber
\mathbb{P}\left(\sup_{0\leq s\leq t}\frac{1}{t^{\frac{1}{\alpha}}}\mid\int_{0}^{s}g(u)dZ^\alpha_0(u)\mid>\varepsilon\right)&=
\mathbb{P}\left(\sup_{0\leq s\leq t}\mid\int_{0}^{s}g(u)dZ^\alpha_0(u)\mid^p>(\varepsilon t^{\frac{1}{\alpha}})^p\right)\\[1ex]
&\leq\frac{3}{(\varepsilon t^{\frac{1}{\alpha}})^p}\mathbb{E}\mid\int_{0}^{t}g(u)dZ^\alpha_0(u)\mid^p.
\end{align}
It follows from \eqref{nier izomt} that
\begin{align}\label{nier. 2}
\frac{3}{(\varepsilon t^{\frac{1}{\alpha}})^p}\mathbb{E}\mid\int_{0}^{t}g(u)dZ^\alpha_0(u)\mid^p\leq\frac{3c_p}{(\varepsilon t^{\frac{1}{\alpha}})^p(p-\alpha)}\int_{0}^{t}\mathbb{E}\mid g(u)\mid^p du.
\end{align}
Since $G$ is Lipschitz, so
\begin{gather}\label{nier. 3}
\mathbb{E}\mid g(u)\mid^p=\mathbb{E}\mid G(R(u-))-G(R(0))\mid^p\leq K\cdot \mathbb{E}\mid R(u-)\mid^p,
\end{gather}
with some constant $K>0$.
By \eqref{nier. 1}, \eqref{nier. 2} and \eqref{nier. 3} we obtain thus
\begin{gather}\label{do dalszego wykorz.}
\mathbb{P}\left(\sup_{0\leq s\leq t}\frac{1}{t^{\frac{1}{\alpha}}}\mid\int_{0}^{s}g(u)dZ^\alpha_0(u)\mid>\varepsilon\right)
\leq \frac{3Kc_p}{(\varepsilon t^{\frac{1}{\alpha}})^p(p-\alpha)}\int_{0}^{t}\mathbb{E}\mid R(u-)\mid^pdu.
\end{gather}
Therefore, for a sequence $\{a_k\}$ we obtain
\begin{align*}
H(k):=\mathbb{P}&\left(\sup_{2^{-k}\leq s\leq 2^{-k+1}}\frac{1}{s^{\frac{1}{\alpha}}}\mid\int_{0}^{s}g(u)dZ^\alpha_0(u)\mid>a_k\right)\\[1ex]
&\leq \mathbb{P}\left(\sup_{2^{-k}\leq s\leq 2^{-k+1}}\frac{1}{(2^{-k+1})^{\frac{1}{\alpha}}}\Big(\frac{2^{-k+1}}{s}\Big)^{\frac{1}{\alpha}}\mid\int_{0}^{s}g(u)dZ^\alpha_0(u)\mid>a_k\right)\\[1ex]
&\leq \mathbb{P}\left(\sup_{0\leq s\leq 2^{-k+1}}2^{\frac{1}{\alpha}}\frac{1}{(2^{-k+1})^{\frac{1}{\alpha}}}\mid\int_{0}^{s}g(u)dZ^\alpha_0(u)\mid>a_k\right)\\[1ex]
&\leq \mathbb{P}\left(\sup_{0\leq s\leq 2^{-k+1}}\frac{1}{(2^{-k+1})^{\frac{1}{\alpha}}}\mid\int_{0}^{s}g(u)dZ^\alpha_0(u)\mid>\frac{a_k}{2^{\frac{1}{\alpha}}}\right), \quad k=0,1,... ,
\end{align*}
and, consequently, by \eqref{do dalszego wykorz.},
\begin{align}\label{przedostatnie szacowanie}
H(k)\leq \frac{3Kc_p}{\Big(\frac{a_k}{2^{\frac{1}{\alpha}}}(2^{-k+1})^{\frac{1}{\alpha}}\Big)^p(p-\alpha)}\int_{0}^{2^{-k+1}}\mathbb{E}\mid R(u-)\mid^pdu,\quad k=0,1,... .
\end{align}
Now we estimate the integral $\int_{0}^{t}\mathbb{E}\mid R(u-)\mid^pdu$ for $t>0$. We can assume that $F$ and $G$ are bounded because we investigate the behaviour of $R$ before it leaves a neighborhood of zero. Then
$$
\mid R(t)\mid^p\leq 2^{p-1}\left(\mid \int_{0}^{t}F(R(s))ds\mid^p+\mid\int_{0}^{t} G(R(s-))dZ^{\alpha}_0(s)\mid^p\right),
$$
and, consequently,
$$
\mathbb{E}\mid R(t)\mid^p\leq 2^{p-1}(ct^p+\mathbb{E}\int_{0}^{t}\mid G(R(s-))\mid^p ds)\leq \tilde{c}t,
$$
with some constants $c,\tilde{c}$. Hence
\begin{gather}\label{szacoewanie na calke z R(u)}
\int_{0}^{t}\mathbb{E}\mid R(u-)\mid^p du=\int_{0}^{t}\mathbb{E}\mid R(u)\mid^p du\leq \tilde{c}\int_{0}^{t}ds=\frac{\tilde{c}}{2}t^2, \quad t>0.
\end{gather}
By \eqref{przedostatnie szacowanie} and \eqref{szacoewanie na calke z R(u)} we obtain finally
\begin{align*}
H(k)&\leq \frac{3Kc_p}{\Big(\frac{a_k}{2^{\frac{1}{\alpha}}}(2^{-k+1})^{\frac{1}{\alpha}}\Big)^p(p-\alpha)}\frac{\tilde{c}}{2}(2^{-k+1})^2\\[1ex]
&=\frac{3Kc_p\tilde{c}2^{\frac{p}{\alpha}-1}}{p-\alpha}\cdot\frac{1}{a_k^p}(2^{-k+1})^{2-\frac{p}{\alpha}},\quad k=0,1,... .
\end{align*}
Taking $a_k=\frac{1}{k}$ and $\delta:=2-p/\alpha>0$ we obtain that
$$
\sum_{k=0}^{+\infty}H_k<+\infty,
$$
and, by the Borel-Cantelli lemma,
$$
\frac{1}{t^{\frac{1}{\alpha}}}\int_{0}^{t}G(R(s-)-G(0))dZ^\alpha_0(s)\underset{t\rightarrow 0}{\longrightarrow} 0,
$$
as required.

\hfill$\square$
\subsection{Analytic HJM condition}
In the  Heath-Jarrow-Morton (HJM)  model of the bond market, see \cite{HeatJarrMor}, the bond prices are written in the form:
\begin{gather}\label{bond prices forward rate}
P(t,T)=e^{-\int_{t}^{T}f(t,s)ds}, \quad 0\leq t\leq T,
\end{gather}
where $f(t,T), 0\leq t\leq T$ is the so-called forward rate process  given by
$$
df(t,T)=\alpha(t,T)dt+\sigma(t,T)dZ(t), \quad f(0,T)=f_0(T), \quad 0\leq t\leq T.
$$
The processes  $\alpha$ and  $\sigma$  are called respectively  drift and volatility of $f$. Moreover, the short rate process equals
\begin{gather}\label{short rate HJM}
R(t)=f(t,t), \quad t\geq 0.
\end{gather}
In the model  introduced by Heath, Jarrow and Morton, $Z$ was a Wiener process. Extensions to models with discontinuous $Z$  were discussed by many authors, see, for instance, \cite{EbJacRaib}, \cite{EbRaib},
\cite{FilTap}, \cite{FilipovicTappeTeichmann}, \cite{BjoDiKabRun},\cite {BjoKabRun}, \cite{Jakubowski-Zabczyk}. In particular,
it was shown that the discounted bond prices, with $T>0$,
$$
\hat{P}(t,T)=e^{-\int_{0}^{t}R(s)ds}e^{-\int_{t}^{T}f(t,s)ds}, \quad t\in[0,T],
$$
are local martingales if and only if
\begin{gather}\label{HJM condition}
A(t,T)=J(\Sigma(t,T)), \quad 0\leq t\leq T,
\end{gather}
for each $T>0$, almost all $t\in[0,T]$ almost $\omega$-surely. Above $A(t,T):=\int_{t}^{T}\alpha(t,s)ds$ and $\Sigma(t,T):=\int_{t}^{T}\sigma(t,s)ds$.

\noindent We transform now  the affine model \eqref{rownanie na P} to the HJM framework and examine the
martingale property via \eqref{HJM condition}. Comparing the exponents in \eqref{rownanie na P} and \eqref{bond prices forward rate} yields
$$
A(T-t)+B(T-t)R(t)=\int_{t}^{T}f(t,s)ds, \quad 0\leq t\leq T,
$$
which determines the forward rate in affine model
\begin{gather}\label{forward rate ATS}
f(t,T)=A^\prime(T-t)+B^\prime(T-t)R(t), \quad 0\leq t\leq T.
\end{gather}
It follows then by setting $T=t$ and \eqref{short rate HJM} that
\begin{gather}\label{warunki na pochodne A i B w zerze}
A^\prime(0)=0, \quad B^\prime(0)=1,
\end{gather}
which, together with \eqref{zerowanie w zerze A i B}, gives preliminary requirements for the functions $A$ and $B$.

Let us define the interval $\bar{I}:=(\bar{a},\bar{b})$, with $0\leq \bar{a}< \bar{b}$, such that
\begin{gather}\label{I dach0}
\bar{x}\in \bar{I}, \qquad G(x)>0  \quad  \text{for} \  x\in\bar{I},
\end{gather}
and such that  $\bar{I}$ is the maximal interval satisfying \eqref{I dach0}.

\begin{prop}\label{tw char. ATS przyp ogolny}
	Assume that, for arbitrary $x\geq 0$,  the non-negative process $R$ is a strong solution of \eqref{rownanie ogolne} with continuous coefficients $F,G$.  Let $A,B$ be twice continuously  differentiable functions satisfying \eqref{zerowanie w zerze A i B} and \eqref{warunki na pochodne A i B w zerze}. Then the affine model \eqref{rownanie na P} has the martingale property if and only if
	\begin{gather}\label{ATS wniosek 1 z HJM ogolny}
	J\Big(G(x)B(v)\Big)=-A^\prime(v)-[B^\prime(v)-1]x+B(v)F(x), \quad x\in\bar{I}, \ v\geq 0.
	\end{gather}
\end{prop}
\vskip2mm

\noindent
{\bf Proof:} We convert the model to the HJM framework. Applying It\^o's formula to \eqref{forward rate ATS} and taking into account
\eqref{rownanie ogolne} we obtain
\begin{align*}
df(t,T)&=-A^{\prime\prime}(T-t)dt+B^\prime(T-t)dR(t)-R(t)B^{\prime\prime}(T-t)dt\\[1ex]
&=\alpha(t,T)dt+\langle \sigma(t,T),dZ(t)\rangle
\end{align*}
where
\begin{align}\label{popostac alpha}
\alpha(t,T)&:=F(R(t)) B^\prime(T-t)-A^{\prime\prime}(T-t)-B^{\prime\prime}(T-t)R(t)\\[1ex]\label{popostac sigma}
\sigma(t,T)&:=B^{\prime}(T-t)G(R(t-)).
\end{align}
Consequently, for $t<T$,
\begin{align*}
A(t&,T)=\int_{t}^{T}\alpha(t,s)ds=\int_{t}^{T}\Big(F(R(t)) B^\prime(s-t)-A^{\prime\prime}(s-t)-B^{\prime\prime}(s-t)R(t)\Big)ds\\[1ex]
&=F(R(t))[B(T-t)-B(0)]-[A^\prime(T-t)-A^\prime(0)]-[B^\prime(T-t)-B^\prime(0)]R(t),\\[1ex]
\Sigma(t,T)&=\int_{t}^{T}\sigma(t,s)ds=\int_{t}^{T}B^{\prime}(s-t)G(R(t-))ds=G(R(t-))[B(T-t)-B(0)].
\end{align*}
Taking into account \eqref{zerowanie w zerze A i B} and \eqref{warunki na pochodne A i B w zerze} we
obtain by \eqref{HJM condition} that the affine model has the martingale property if and only if
\begin{gather}\label{basic equation from HJM}
J\Big(G(R(t-))B(T-t)\Big)=-A^\prime(T-t)-[B^\prime(T-t)-1]R(t)+B(T-t)F(R(t))
\end{gather}
for each $T>0$, $\mathbb{P}$-almost surely, for almost all $t\in[0,T]$.

Now we prove that \eqref{basic equation from HJM} is equivalent to \eqref{ATS wniosek 1 z HJM ogolny}.
Since, for almost all $t\geq 0$, $R(t-)=R(t)$ one can replace $R(t-)$ by $R(t)$ in \eqref{basic equation from HJM}. So, it is
clear that \eqref{ATS wniosek 1 z HJM ogolny} is sufficient for \eqref{basic equation from HJM}. Now we show necessity.
Assume to the contrary that for some $\bar{x}>0$ and $\bar{v}>0$
\begin{gather*}
J(G(\bar{x})B(\bar{v}))>-A^\prime(\bar{v})-[B^\prime(\bar{v})-1]\bar{x}+B(\bar{v})F(\bar{x}).
\end{gather*}
Then there exists $\delta>0$ such that
\begin{gather*}
J(G(x)B(v))>-A^\prime(v)-[B^\prime(v)-1]x+B(v)F(x),
\end{gather*}
for $x\in(\bar{x}-\delta,\bar{x}+\delta)$ and $v\in(\bar{v}-\delta,\bar{v}+\delta)$. Let us consider
the solution $R$ of \eqref{rownanie ogolne} starting from $\bar{x}$ and let us define
$$
\tau:=\inf\{t\geq 0: \mid R(t)-\bar{x}\mid>\delta\}.
$$
Then for $t\in(0,\tau)$  and $T$ such that $T-t\in(\bar{v}-\delta,\bar{v}+\delta)$ 
$$
J(G(R(t-))B(T-t))>-A^\prime(T-t)-[B^\prime(T-t)-1]R(t)+B(T-t)F(R(t)),
$$
which is a contradiction.\hfill$\square$

\section{Proofs of the main results}
\subsection{Proof of Part (I) of Theorem \ref{tw glowne necessity}}

The proof is rather involved and therefore is divided into several steps. The general idea is as follows. First one proves the affine formula for $F$ on some subinterval of $[0,+\infty)$. Then one establishes that the process $Z$ is an $\alpha$-stable martingale with $\alpha \in (1,2]$.  For this one proves first that on an open subinterval of $(0, +\infty)$,  $G(x) = c(x+B)^{\gamma}$, with some constants $c, B, \gamma$, and then one deduces that $J(\lambda)= c_\alpha\lambda^{\alpha}$ first locally and then globally on the whole $[0,+\infty )$. In the final step one shows that $G(0)=0$ and thus that $B=0$.
\vskip2mm

Let  $\bar{I}:=(\bar{a},\bar{b})$,  $0\leq \bar{a}< \bar{b}$, be the maximal interval for which
\begin{gather}\label{I dach}
\bar{x}\in \bar{I}, \qquad G(x)>0  \quad  \text{for} \  x\in\bar{I}.
\end{gather}
If the affine model with short rate $R$ and functions $A$, $B$ has the martingale property then it follows from
Theorem \ref{tw char. ATS przyp ogolny} that
\begin{gather}\label{ATS wniosek 1 z HJM ogolny w dowodzie}
J\Big(G(x)B(v)\Big)=-A^\prime(v)-[B^\prime(v)-1]x+B(v)F(x), \quad v\geq 0, \ x\in\bar{I}.
\end{gather}
\vskip2mm
\noindent
\underline{{\bf Step 1:}} We prove the linear form of $F$ on the interval $\bar{I}$.
\vskip2mm

Differentiation of \eqref{ATS wniosek 1 z HJM ogolny w dowodzie} yields
$$
J^\prime(G(x)B(v))G(x)B^\prime(v)=-A^{\prime\prime}(v)-B^{\prime\prime}(v)x+B^\prime(v)F(x), \quad v\geq 0, \ x\in\bar{I}.
$$
Putting $v=0$ yields
$$
J^\prime(0+)G(x)=-A^{\prime\prime}(0)-B^{\prime\prime}(0)x+F(x), \quad \ x\in\bar{I}.
$$
Since $J^\prime(0+)=0$ we obtain the formula for $F$
\begin{gather}\label{dryf liniowy w dowodzie}
F(x)=A^{\prime\prime}(0)+B^{\prime\prime}(0)x:=ax+b, \quad \ x\in\bar{I}.
\end{gather}

\noindent  To show that\,\, $b\geq 0 ,$\,\, assume, by contradiction, that $b< 0$ and consider solution $R$ of the equation starting from $0$. Since $R$ is non-negative, we have:
\[
R(t) = b\int_0^t e^{a(t-s)}ds + \int_0^t e^{a(t-s)}G(R(s-))dZ(s) \geq 0,\,\,t\geq 0.
\]
Hence
\[
-|b| \int_0^t e^{-as}ds + \int_0^t e^{-as}G(R(s-))dZ(s)\geq 0, \quad t\geq 0,
\]
or equivalently
\[
\int_0^t e^{-as}G(R(s-))dZ(s)\geq |b| \int_0^t e^{-as}ds\quad t\geq 0.
\]
Since the above  stochastic integral is a local non-negative martingale starting from $0$ it must be identicaly $0$. Thus the process
$$
R(t)= b\int_0^t e^{a(t-s)}ds,\,\,t\in [0,T]
$$
is strictly negative and we have contradiction.
\vskip2mm

\noindent
\underline{{\bf Step 2:}} We prove that $G$ satisfies
\begin{gather}\label{G prime przez G}
\frac{G^\prime(x)}{G(x)}=\frac{\bar{A}}{\bar{B}+x}, \quad x\in\bar{I}, \quad x\neq \bar{B},
\end{gather}
with some constants $\bar{A},\bar{B}\in\mathbb{R}, \bar{A}\neq0$.
\vskip2mm

Since \,\,$F(x) = ax +b , x\in\bar{I}$,
\eqref{ATS wniosek 1 z HJM ogolny w dowodzie} yields
\begin{gather}\label{ATS wniosek 1 z HJM ogolny w dowodzie F liniowe}
J(G(x)B(v))=-A^\prime(v)+B(v)b+x[B(v)a+1-B^\prime(v)],\quad v\geq 0, \ x\in\bar{I}.
\end{gather}
Differentiation of \eqref{ATS wniosek 1 z HJM ogolny w dowodzie F liniowe} over $x$ and $v$ leads to
\begin{align}\label{Row 1}
J^\prime(G(x)B(v))G^\prime(x)B(v)&=B(v)a+1-B^\prime(v),\quad v\geq 0, \ x\in\bar{I},\\[1ex]\label{Row 2}
J^\prime(G(x)B(v))G(x)B^\prime(v)&=-A^{\prime\prime}(v)+B^\prime(v)b+x[B^\prime(v)a-B^{\prime\prime}(v)],\quad v\geq 0, \ x\in\bar{I}.
\end{align}
Since $B$ is continuously differentiable and $B^\prime(0)=1$, we can find an $\varepsilon>0$ such that
$$
B(v)>0, \quad B^\prime(v)>0, \quad v\in(0,\varepsilon).
$$
Let us assume that the right side of \eqref{Row 1} equals zero for $v\in(0,\varepsilon)$. Then $B$ solves
$$
B^\prime(v)=aB(v)+1, \quad B(0)=0,\quad v\in(0,\varepsilon),
$$
so is, on the interval $(0,\varepsilon)$, of the form $B(v)=\frac{1}{a}(e^{av}-1)$ if $a\neq 0$ or $B(v)=v$ if $a=0$.
Since the left side of \eqref{Row 1} equals zero and $B(v)> 0$ for $v\in(0,\varepsilon)$ and $G^\prime(\bar{x})\neq 0$, we obtain
$$
J^\prime(G(\bar{x})B(v))=0, \quad v\in(0,\varepsilon).
$$
Hence $J^\prime$ disappears on some interval and consequently must disappear on $[0,+\infty)$ as an analytic function.
Since $J(0)=0$, this implies that $J(\lambda)=0$ for $\lambda\in [0,+\infty)$, which is impossible.
It follows thus that the right side of \eqref{Row 1} is different from zero for some $\bar{v}\in(0,\varepsilon)$. This implies that
$$
B(\bar{v})\neq 0, \quad G^\prime(x)\neq 0,\quad J^\prime(G(x)B(\bar{v}))\neq 0, \quad x\in\bar{I}.
$$
Hence, by \eqref{Row 1},
$$
J^\prime(G(x)B(\bar{v}))=\frac{B(\bar{v})a+1-B^\prime(\bar{v})}{G^\prime(x)B(\bar{v})}, \quad x\in \bar{I}.
$$
Putting this into \eqref{Row 2} with $v=\bar{v}$ yields
$$
\frac{B(\bar{v})a+1-B^\prime(\bar{v})}{G^\prime(x)B(\bar{v})} \cdot G(x)B^\prime(\bar{v})=
-A^{\prime\prime}(\bar{v})+B^\prime(\bar{v})b+x[B^\prime(\bar{v})a-B^{\prime\prime}(\bar{v})], \quad x\in \bar{I},
$$
and, consequently,
\begin{gather}\label{roownanie na G}
\frac{G(x)}{G^\prime(x)}=\frac{-A^{\prime\prime}(\bar{v})+B^\prime(\bar{v})b+x[B^\prime(\bar{v})a-B^{\prime\prime}(\bar{v})]}{B(\bar{v})a+1-B^\prime(\bar{v})}\cdot \frac{B(\bar{v})}{B^\prime(\bar{v})}, \quad x\in \bar{I}.
\end{gather}
Hence the quotient $G(x)/G^\prime(x)$ is a linear function of $x$. If $B^\prime(\bar{v})a-B^{\prime\prime}(\bar{v})=0$ then
$$
\frac{G(x)}{G^\prime(x)}=c, \quad x\in\bar{I},
$$
with $c=G(\bar{x})/G^\prime(\bar{x})\neq 0$. Thus in this case $G(x)=ke^{\frac{x}{c}}, x\in \bar{I}$ with $k>0$. Using this form of $G$ in
\eqref{ATS wniosek 1 z HJM ogolny w dowodzie F liniowe} gives
$$
J(ke^{\frac{x}{c}}B(v))=-A^\prime(v)+B(v)b+x[B(v)a+1-B^\prime(v)],\quad v\geq 0, \ x\in\bar{I},
$$
and allows us to determine $J$ by
$$
J(\lambda)=-A^\prime(v)+B(v)b+\ln\lambda \frac{c}{B(v)k}[B(v)a+1-B^\prime(v)], \quad \lambda>0, v\geq0.
$$
The right side above must be independent of $v$ and  $J(0+)=\pm\infty$ or $J(0)\equiv 0$. Both situations are not possible, so
we conclude  that $B^\prime(\bar{v})a-B^{\prime\prime}(\bar{v})\neq 0$ and thus \eqref{G prime przez G} holds.

\vskip2mm

\noindent
\underline{{\bf Step 3:}} We prove that $Z$ is a stable martingale with index $\alpha\in(0,2]$ by examining
\eqref{G prime przez G} with non-negative and negative $\bar{B}$.

\vskip2mm
$(a)$ If $\bar{B}\geq 0$ then \eqref{G prime przez G} can be written in the form
$$
\frac{d}{dx}\ln(G(x))=\bar{A}\frac{d}{dx}\ln(\bar{B}+x), \quad x\in\bar{I},\quad x\neq \bar{B},
$$
which yields
$$
\ln G(x)-\ln(\bar{B}+x)^{\bar{A}}=k, \quad x\in\bar{I},\quad x\neq \bar{B},
$$
for some $k\in\mathbb{R}$. Consequently,
\begin{gather}\label{G na poddprzedfziale}
G(x)=K(\bar{B}+x)^{\bar{A}}, \quad x\in\bar{I},\quad x\neq \bar{B},
\end{gather}
with $K:=e^{k}$. Now we put \eqref{G na poddprzedfziale} into \eqref{ATS wniosek 1 z HJM ogolny w dowodzie F liniowe}. This yields
$$
J\Big(K(\bar{B}+x)^{\bar{A}} B(v)\Big)=-A^\prime(v)+B(v)b+x[B(v)a+1-B^\prime(v)],\quad v\geq 0, \ x\in\bar{I}, \ x\neq\bar{B}.
$$
We fix $v=\tilde{v}$ such that $B(\tilde{v})\neq 0$ and introduce $z:=K(\bar{B}+x)^{\bar{A}} B(\tilde{v})$. Then $x=(K B(\tilde{v}))^{-\frac{1}{\bar{A}}}z^{\frac{1}{\bar{A}}}-\bar{B}$ and consequently
$$
J(z)=k_1z^{\frac{1}{\bar{A}}}+k_2, \quad z\in \bar{J}:=(K(\bar{B}+\bar{a})^{\bar{A}} B(\tilde{v}),K(\bar{B}+\bar{b})^{\bar{A}} B(\tilde{v}))
$$
with some constants $k_1,k_2$. Since $J(0)=0$ and $J^\prime(0+)<+\infty$ we obtain that
$$
J(z)=k_1z^{\frac{1}{\bar{A}}},\quad z\in \bar{J},
$$
and $\alpha:=\frac{1}{\bar{A}}>1$. In fact, $\alpha$ can not be greater than $2$, which we show below.
By \eqref{J wzor}, for $Z$ without negative jumps,
\begin{align*}
J(z)&=-az+\frac{1}{2}qz^2+\int_{0}^{+\infty}\left(e^{-zy}-1+zy\mathbf{1}_{(-1,1)}(y)\right)\nu(dy)\\[1ex]
&=-a z+\frac{1}{2}qz^2+\int_{0}^{1}\left(e^{-zy}-1+zy\right)\nu(dy)+\int_{1}^{+\infty}\left(e^{-zy}-1\right)\nu(dy)\\[1ex]
&=-a z+\frac{1}{2}qz^2+z^2\int_{0}^{1}\frac{e^{-zy}-1+zy}{(zy)^2}y^2\nu(dy)+\int_{1}^{+\infty}\left(e^{-zy}-1\right)\nu(dy),\quad z\geq 0.
\end{align*}
Since the function
$$
x\rightarrow \frac{e^{-x}-1+x}{x^2},\quad  x\geq 0,
$$
is bounded, the measure $y^2\nu(dy)$ is finite on $[0,1]$ and
$$
\int_{1}^{+\infty}\left(e^{-zy}-1\right)\nu(dy)\leq z \int_{1}^{+\infty}y\nu(dy),
$$
we see actually that $J(z)\leq az^2+bz+c$ for some positive constants $a,b,c$. Since $J$ is analytic on $[0,+\infty)$ we obtain finally
$$
J(z)=c_{\alpha}z^{\alpha},\quad z\in [0,+\infty),
$$
where $c_\alpha=k_1$. Therefore $Z$ is a stable process with index $\alpha\in(1,2]$.

$(b)$ If $\bar{B}<0$ and $-\bar{B}\leq \bar{x}$ then we can examine \eqref{G prime przez G} on the set $\{x:x\in\bar{I}, x>-\bar{B}\}$ in the same way as in  $(a)$.
If $\bar{B}<0$ and $-\bar{B}> \bar{x}$ then we examine \eqref{G prime przez G} on the set $\{x:x\in(\bar{x},-\bar{B}\wedge \bar{b})\}$ and write it in the form
\begin{gather*}
\frac{G^\prime(x)}{G(x)}=\frac{-\bar{A}}{-\bar{B}-x}, \quad x\in(\bar{x},-\bar{B}\wedge \bar{b}).
\end{gather*}
It follows that
$$
\frac{d}{dx}\ln(G(x))=\bar{A}\frac{d}{dx}\ln(-\bar{B}-x), \quad x\in(\bar{x},-\bar{B}\wedge \bar{b})
$$
and, consequently
\begin{gather}\label{G na podprzedfziale}
G(x)=K(-\bar{B}-x)^{\bar{A}}, \quad x\in(\bar{x},-\bar{B}\wedge \bar{b}),
\end{gather}
with some $K>0$.
Putting this formula into \eqref{ATS wniosek 1 z HJM ogolny w dowodzie F liniowe} yields
$$
J\Big(K(-\bar{B}-x)^{\bar{A}} B(v)\Big)=-A^\prime(v)+B(v)b+x[B(v)a+1-B^\prime(v)],\quad v\geq 0, \quad x\in(\bar{x},-\bar{B}\wedge \bar{b}).
$$
Now, like in $(a)$ one shows that $J(z)=c_\alpha z^{\alpha}, z\in[0,+\infty)$ with $\alpha:=1/\bar{A}\in(1,2]$.

\vskip2mm

\noindent
\underline{{\bf Step 4}}: We prove that $\bar{I}=(0,+\infty)$ and that
\begin{gather}\label{wzor na G final}
G(x)=c^{\frac{1}{\alpha}}x^{\frac{1}{\alpha}}, \qquad c>0, \ x\in[0,+\infty).
\end{gather}
Since we know from Step 3 that $J(z)=c_\alpha z^{\alpha}, z\in[0,+\infty)$, $\alpha\in(1,2]$, it follows from
\eqref{Row 1} that
$$
\alpha c_\alpha G^{\alpha-1}(x)G^\prime(x)B^\alpha(v)=B(v)a+1-B^\prime(v), \quad x\in\bar{I}, \quad v\geq 0.
$$
We can find $\tilde{v}>0$ such that $B(\tilde{v})\neq 0$. Then
\begin{gather}\label{krokodyl zielony}
\alpha c_{\alpha} G^{\alpha-1}(x)G^\prime(x)=M, \quad x\in{\bar{I}},
\end{gather}
with $M:=(B(\tilde{v})a+1-B^\prime(\tilde{v}))/B^\alpha(\tilde{v})$. Now we show that $\bar{I}=(0,+\infty)$.
Assume that $\tilde{a}>0$. Since, by definition, $\lim_{x\downarrow \bar{a}}G(x)=0$, we see from \eqref{krokodyl zielony} that
$\lim_{x\downarrow \bar{a}}G^\prime(x)=\pm \infty$, which contradicts the differentiability of $G$ on $(0,+\infty)$.
Similarly one can exclude the case $\tilde{b}<+\infty$. Solving \eqref{krokodyl zielony} we obtain
$$
G(x)=\Big(G(\bar{x})-\frac{M}{c_{\alpha}}\bar{x}+\frac{M}{c_{\alpha}}x\Big)^\frac{1}{\alpha}:=(m_1+m_2x)^{\frac{1}{\alpha}}, \quad x\in(0,+\infty),
$$
with $m_1\geq 0, m_2>0$. If $m_1>0$ then $G$ is Lipschitz at zero and by Propositon \ref{tw G(0)=0}, $G(0)=0$ which is a contradiction.
Hence \eqref{wzor na G final} follows with $c:=m_2$.

\hfill $\square$

\subsection{{Proof of Part (II) of Theorem \ref{tw glowne necessity}}}
By elementary arguments,  positivity of the solutions to the stochastic equation \eqref{rownanie ogolne} with $G(x)\equiv\sigma$ implies that $Z$ has no Wiener part and can have only positive jumps. Repeating the arguments from the proof of Part (I) one can show that  
that $F(x) = ax + b , x\geq 0$, and $b\geq 0$. We will establish now that
\begin{equation}
  \int_0^{+\infty}y\nu (dy) < +\infty, \quad b\geq \sigma \int_0^{+\infty}y\nu (dy).
\end{equation}
Let $\tilde \pi$ \,\,be the compensated jump measure corresponding to the martingale $Z$. Then for $\epsilon >0$,
\begin{align}
Z(t) = &\int_0^t \int_0^{+\infty} y \tilde \pi (ds, dy)=\int_0^t \int_0^{\epsilon} y \tilde \pi (ds, dy) +\int_0^t \int_{\epsilon}^{+\infty} y \tilde \pi (ds, dy)\\
= & Z_{\epsilon} (t) + P_{\epsilon}(t) - t \int_{\epsilon}^{+\infty} y\nu(dy).
\end{align}
Here $Z_{\epsilon}$ is a L\'evy martingale with positive jumps bounded by $\epsilon$, $P_{\epsilon}$ is a compound Poisson process with the L\'evy measure $\nu$ restricted to the interval $[\epsilon , +\infty )$. For the solution $R$ of the stochastic equation, starting from $0$, we have for all $t\in [0, T]$:
\begin{align}\label{rozlozone rownanie}\nonumber
e^{-at}R(t) =& b \int_0^t e^{-as}ds + \sigma \int_0^t e^{-as}dZ(s)\\
= &(b- \sigma \int_{\epsilon}^{+\infty}y\nu(dy))\int_{0}^{t}e^{-as}ds + \sigma \int_{0}^{t}e^{-as}dZ_{\epsilon}(s) +\int_{0}^{t}e^{-as}dP_{\epsilon}(s)\geq 0.
\end{align}
If $\int_{0}^{+\infty} y \nu(dy) = +\infty$, then by taking $\epsilon$ close to $0$, the number $(b- \sigma \int_{\epsilon}^{+\infty}y\nu(dy))$ can be made arbitrary small negative. The stochastic integrals with respect to $Z_{\epsilon},\, P_{\epsilon}$ are independent processes. The former one can be made, with positive probability, uniformly smaller on $[0, T]$, than given in advance number and the latter one , is $0$ on $[0,T]$
also with positive probability. Thus $e^{-at}R(t)$ is negative for some $t\in [0,T]$, with positive probability, which is a contradiction.
Now we show that $b\geq \sigma\int_{0}^{+\infty}y\nu(dy)$. In the opposite case we have that the difference
$$
b- \sigma\int_{\epsilon}^{+\infty}y\nu(dy),
$$
is negative for sufficiently small $\epsilon>0$ and decreases as $\epsilon\downarrow 0$. It follows from the Markov inequality that for any $\gamma>0$ and $t>0$
$$
\mathbb{P}(\sigma\int_{0}^{t}e^{-as}dZ_{\epsilon}(s)>\gamma)\leq
\frac{\sigma^2\mathbb{E}(\int_{0}^{t}e^{-as}dZ_{\epsilon}(s))^2}{\gamma^2}=
\frac{\sigma^2\int_{0}^{t}\int_{0}^{\epsilon}e^{-2as}y^2ds\nu(dy)}{\gamma^2}\underset{\epsilon\rightarrow 0}{\rightarrow} 0,
$$
and consequently
$$
\mathbb{P}(\sigma\int_{0}^{t}e^{-as}dZ_{\epsilon}(s)\leq \gamma)\underset{\epsilon\rightarrow 0}{\rightarrow} 1.
$$
Since the integral over $P_\varepsilon$ disappears with positive probability, we have by \eqref{rozlozone rownanie} that $R(t)<0$ which is a contradiction.\hfill$\square$

\subsection{Proof of Part (I) of Theorem \ref{tw sufficiency}}
It was shown in \cite{FuLi} that equation \eqref{rownanie 1002} actually has a unique non-negative strong solution.
Now we use Theorem \ref{tw char. ATS przyp ogolny} with $J(\lambda)=c_{\alpha}\lambda^{\alpha}$, $F(x)=ax+b$ and
$G(x)=c^{\frac{1}{\alpha}}x^{\frac{1}{\alpha}}$. Then \eqref{ATS wniosek 1 z HJM ogolny} boils down to
$$
c_{\alpha}\Big(c^{\frac{1}{\alpha}}x^{\frac{1}{\alpha}}B(v)\Big)^\alpha=-A^\prime(v)-[B^\prime(v)-1]x+B(v)[ax+b], \quad x\geq 0, \quad v\geq 0.
$$
Consequently,
\begin{gather}\label{do analizy}
c_{\alpha}cxB^\alpha(v)=(aB(v)-B^\prime(v)+1)x+bB(v)-A^\prime(v), \quad x\geq 0, \quad v\geq 0.
\end{gather}
Putting $x=0$ yields
$$
bB(v)-A^\prime(v)=0, \quad v\geq 0,
$$
which is the required formula for $A$. It follows from \eqref{do analizy} that
$$
c_{\alpha}cB^\alpha(v)=aB(v)-B^\prime(v)+1, \quad v\geq 0,
$$
which yields the equation for $B$.
\hfill$\square$
\subsection{Proof of Part (II) of Theorem \ref{tw sufficiency}}
Note that functions $A, B$ should satisfy, for all $x\geq 0,\,v\geq 0$,\, the equation
\begin{align*}
J(\sigma B(v))=& -A'(v) -(B'(v) - 1)x + B(v) (ax+b)\\
=& x (B(v)a - B'(v) +1 ) + B(v)b - A'(v).
\end{align*}
Consequently
\begin{align*}
 B'(v)= &B(v)a +1,\quad B(0)=0\\
A'(v)=& B(v)b -J(\sigma B(v)), \quad A(0)=0.
\end{align*}
It remains to show that $A$ is an increasing function that is, that is  $B(v)b -J(\sigma B(v))\geq 0.$

\noindent However,
\begin{align}
J(\lambda) = &\int_{0}^{+\infty} (e^{-\lambda y} -1)\nu(dy) + \lambda \int_{0}^{+\infty}y\nu(dy)\\
=:& J_{0}(\lambda) + \lambda \int_{0}^{+\infty}y\nu(dy).
\end{align}
It is clear that
\[
B(v)(b - \sigma \int_{0}^{+\infty}y\nu(dy)) - J_{0}(\sigma B(v)) \geq 0 ,
\]
and the result follows.
\hfill$\square$

\end{document}